\documentclass[10pt]{article}
\font\smallit=cmti10

\usepackage{tabularx}
\usepackage{amssymb,amsmath,amsthm,enumitem} 
\usepackage{subcaption,booktabs}

\makeatletter

\renewcommand\section{\@startsection {section}{1}{\z@}
{-30pt \@plus -1ex \@minus -.2ex}
{2.3ex \@plus.2ex}
{\normalfont\normalsize\bfseries\boldmath}}

\renewcommand\subsection{\@startsection{subsection}{2}{\z@}
{-3.25ex\@plus -1ex \@minus -.2ex}
{1.5ex \@plus .2ex}
{\normalfont\normalsize\bfseries\boldmath}}

\renewcommand{\@seccntformat}[1]{\csname the#1\endcsname. }

\makeatother

\newtheorem{theorem}{Theorem}
\newtheorem{lemma}{Lemma}

\newtheorem{conjecture}{Conjecture}

\theoremstyle{definition}

\newtheorem{remark}{Remark}

\def\R{\mathcal{R}}
\def\M{\mathcal{M}}

\def\ZZ{\mathbb{Z}}

\begin{document}

\begin{center}
\uppercase{$\mathcal{GM}$-rule and its applications to impartial games}
% \uppercase{\bf M-$\ell$-Rule and Remoteness Functions
% of Exact Slow $\MakeLowercase{k}$-NIM with $\MakeLowercase{k}+1$ Piles}

\vskip 20pt
{\bf Vladimir Gurvich}\\
{\smallit National Research University Higher School of Economics (HSE), Moscow, Russia}\\
{\tt vgurvich@hse.ru}, {\tt vladimir.gurvich@gmail.com}\\
\vskip 10pt
% {\bf Vladislav Maximchuk}\\
% {\smallit National Research University Higher School of Economics (HSE), Moscow, Russia}\\
% {\tt vladislavmaximchuk3495@gmail.com}\\
% \vskip 10pt
% {\bf Georgy Miheenkov}\\
% {\smallit National Research University Higher School of Economics (HSE), Moscow, Russia}\\
% {\tt georg2002g@mail.ru}\\
% \vskip 10pt
{\bf Mariya Naumova}\\
{\smallit Rutgers Business School, Rutgers University, Piscataway, NJ, United States}\\
{\tt mnaumova@business.rutgers.edu}\\
% \vskip 10pt
% {\bf Michael Vyalyi}\\
% {\smallit National Research University Higher School of Economics, Moscow, Russia;\\
%  Institute of Physics and Technology, Dolgorpudnyi, Russia;\\
% Federal Research Center 
% ``Computer Science and Control'' 
% of the Russian Academy of Sciences, Moscow, Russia}\\
% {\tt vyalyi@gmail.com}
\end{center}
\vskip 20pt
\centerline{\smallit Received: , Revised: , Accepted: , Published: } 
% We will fill in the dates
\vskip 30pt

\centerline{\bf Abstract}
\noindent
Given integer $n \geq 1, \ell \geq 2$, 
and vector $x = (x_1, \ldots, x_n)$  
that has an entry which is a multiple of $\ell$ and 
such that $x_1 \leq \ldots \leq x_n$,   
the GM-rule is defined as follows: 
Keep the rightmost 
minimal entry  $x_i$ of $x$, which is a multiple of $\ell$ 
% with the largest index 
% (it is called the {\em pivot} of $x$) 
and reduce the remaining $n-1$ entries of $x$ by~1. 
We will call such % (unique) 
$i$  the {\em pivot} 
and $x_i$  the {\em pivotal entry}. 
The GM-rule respects monotonicity of the entries. 
It uniquely determines a GM-move $x^0 \to x^1$  and 
an infinite GM-sequence $S$ that consists 
of successive GM-moves $x = x^0 \to x^1 \to \ldots \to x^j \to \ldots$     
% $ \max(x_i \mid i \in [n]) - \min(x_i \mid i \in [n])$, 
If $range(x) = x_n - x_1 \leq \ell$ then for all $j \geq 0$:
% we have:  
\begin{itemize}
\item[(i)]  $range(x^j) \leq \ell$;    
\item[(ii)] the pivot of  $x^{j + \ell}$  is one less than 
the pivot of  $x^j$, assuming that $1 - 1 = 0 = n$. 
\item[(iii)] % GM-sequence has the period  $n \ell$; more precisely  
$x_i^j  - x_i^{j + n \ell} = (n-1) \ell$   
for all  $i = 1,\ldots,n$. %  and $j \geq 0$; 
\end{itemize}
Due to  (iii), we compute  $x^j$  
in time linear in  
$n, \ell, \log(j)$, and $\sum^n_{i=1}\log(|x_i|+1)$. 
% $\sum_{i \in [n]}(\log(|x_i| + 1)$

For $\ell = 2$  a slighty modified version 
of the GM-rule was recently introduced by 
Gurvich, Martynov, Maximchuk, and Vyalyi,  
``On Remoteness Functions of 
\newline 
Exact Slow $k$-NIM with $k+1$ Piles",  
https://arxiv.org/abs/2304.06498 (2023), 
where applications to impartial games were considered. 

% \newline 
% \noindent{\bf AMS subjects:  91A05, 91A46, 91A68} 
% \pagestyle{myheadings}
% \markright{\smalltt  INTEGERS: 23 (2023)\hfill}
% \thispagestyle{empty}
% \baselineskip=12.875pt
% \vskip 30pt
%%Your paper starts here with your first section.

\section{Introduction} 
\label{s0}
Fix integer $n \geq 1, \ell \geq 2$, 
and consider an integer vector $x = (x_1, \ldots, x_n)$. Note that $x$ may have negative entries.  
We assume that $x$ is defined 
up to a permutation of its enties. 
So without loss of generality (wlog) 
we will require monotonicity:  
$x_1 \leq \ldots \leq x_n$. 
Assume also that  $x$  has an entry which is a multiple of $\ell$. 
Choose a minimal entry of $x$  
which is a multiple of $\ell$. 
If there are several such entries,  
take one with the largest index. 
This index  $i \in [n] = \{1,\ldots,n\}$  
will be called the {\em pivot} 
and $x_i$ the {\em pivotal entry} of $x$. 
% Obviously, the pivot is unique.  

\medskip 

The {\em GM-rule} is defined as follows: 
keep the pivotal entry unchanged 
and reduce the remaining $n-1$ entries by 1.
This transformation is called the {\em GM-move} 
(from $x = x^0$)  and denoted  $x^0 \to x^1$.   
The (unique) infinite sequence of successive GM-moves \;\;
% from~$x^0$  \; 
$S = S(x^0) = (x^0 \to x^1 \to \ldots \to x^j \to \ldots)$ 
is called the {\em GM-sequence}.
% (from $x*0$)  and denoted by  $S = S(x^0)$.
Obviously, GM-moves respect monotonicity, 
as well as the existence of an entry which is a multiple of~$\ell$. 
Indeed, it may disappear 
only after a replacement (a new pivot) appears.

\begin{remark}
If $x$ has no entry which is a multiple of $\ell$, 
we can introduce the GM-move in different ways: 
it can keep $x_n$ (as in \cite{GMMV23})   
or it can reduce each entry by 1. 
In any case, such vector  $x$ is irreturnable, 
since, obviously, no GM-move can enter it. 
Furthermore, after at most $\ell-1$  GM-moves 
from $x$  the GM-sequence $S(x)$  comes 
to a vector $x'$  
that has an enty which is a multiple of $\ell$.  
In the sequel, we assume everywhere, 
except Section \ref{s-games}, 
that $x$  has this property. 
% an entry multiple of $\ell$. 
\end{remark}
\medskip
   
Given a vector $x = (x_i \mid i \in [n])$, 
define its {\em range} by formula 
$range(x) = x_n - x_1$.
%% \max(x_i \mid i \in [n]) - \min(x_i \mid i \in [n])$, 
%% where $[n] = \{1, \ldots, n\}$.  
A GM-move may change the range by at most 1:  
the range is reduced by 1 if the pivot is $1$,   
increased by 1 if the pivot is  $n$, and 
remains the same if pivot is $i$  such that $1 < i < n$. 
Assuming that $range(x^0) \leq \ell$,  
in Section \ref{s-main},  
we will prove the following properties of a 
GM-sequence $S = S(x^0)$. For each $j \geq 0$ we have:    
% =(x^0 \to x^1 \to \ldots \to x^j \to \ldots)$: 

\begin{itemize}
\item[(i)]  $range(x^j) \leq \ell$; %  for every $j$  in $S$;   
\item[(ii)] the sequence of pivots has the period $\ell$; 
more precisely, the pivot of  $x^{j+\ell}$  is one less than 
the pivot of $x^j$, 
assuming standardly that $n$ is one less than 1. 
\item[(iii)] % GM-sequence 
$S$ has the period $n \ell$; 
more precisely,  $x_i^{j + n \ell} = x_i^j - (n-1) \ell$   
for all  $i \in [n]$.  % and $j \geq 0$. 
\end{itemize}

In Section \ref{s-main}, based on (iii), 
we compute $x^j$ in time linear in  
$n, \ell, \log(j)$, and 
$\omega(x) = \sum_{i \in [n]}(\log(|x_i| + 1)$.
% The algorithm is based on property (iii). 

In subsection \ref{ss-racing} we prove the following statement: 

\begin{itemize}
\item[(iv)] 
If $range(x^0) > \ell$ then 
there exists an integer $N = N(S)$,   
linearly bounded in  $n, \, \ell$ and $x$, 
% $\sum_{i \in [n]}(\log(|x_i| + 1)$ 
such that $range(x^j) \leq \ell$ if and only if  $j > N$. 
Furthermore, $N$  will be computed in subsection \ref{ss-racing} 
in time linear in $n, \, \ell$, and  $\omega(x)$.
\end{itemize}

We can apply (i, ii, iii) just replacing $x^0$ by $x^{N+1}$.   

Thus, an SG-sequence $S = S(x^0)$ partitions  $\ZZ_+$ into  $S_0 = \{0, 1, \ldots, N\}$ and $S_\infty =\{N+1,N+2,\ldots\}$. 

Then, (i) and (iv) imply that set $S_\infty$ is absorbing. 
Note that $S_0$  may be empty;  
in this case we set $N = 0$, by convention. 

\medskip 

For $\ell = 2$, a slightly different version of the GM-rule 
was introduced in \cite{GMMV23} and called the {\em M-rule}. 
We view the GM-rule as a generalization of it.

In \cite{GMMV23}, an impartial game 
``Exact Slow $k$-NIM with $n$ Piles"  
was considered for $n = k+1$. 
In this game positions $x$  are non-negative,  
while  $\M(x)$  is the number of moves in GM-sequence  $S(x)$ 
until at least two of $n$ entries become  0. 
(No entry can become negative. 
Indeed, as soon as one entry is reduced to $0$ 
it becomes a pivot and is not reduced further;  
as soon as two entries become 0, the game is over.) 

In \cite{GMMV23} it was proven that $\R(x) = \M(x)$, 
where  $\R(x)$  is the so-called 
{\em remoteness function}  \cite{Smi66} 
of the game in position $x$. 
The value $\R(x)$  equals the (unique) length 
of an optimal play, 
assuming that the winner wants to win as soon as possible, 
while the loser wants to resist as long as possible. 
Furthermore, $\R(x)$ is even 
if and only if $x$ is a P-position. 
See precise definitions and more detail  
in \cite{GMMV23} or in Section \ref{s-games}.  
In \cite{GMMV23} it was proven that 
$\R(x) = \M(x)$ and the P-positions  
were computed in time linear in $n, \ell$, and $\omega(x)$. 
% $\sum_{i \in [n]}(\log(|x_i| + 1)$. 
Somewhat surprisingly,  
explicit formulas for the P-positions 
are known only for $n \leq 4$ 
(see \cite{GHHC20} and \cite[Appendix]{GMMV23})  
and they are pretty complicated. 
In this paper we construct a new algorithm 
computing  $\M(x)$ in time linear in $n, \ell$, and  $\omega(x)$.

\section{Main results} 
\label{s-main}
\subsection{Range $\ell$ lemma}
First, let us prove that the set 
$\{x \mid range(x) \leq \ell\}$  is absorbing. 

\begin{lemma} 
\label{l2} 
For each GM-sequence  
that begins in $x$ and ends in $x'$  we have:  

\smallskip 

(i) if $range(x) \leq \ell$  then  $range(x') \leq \ell$;

\smallskip 

(ii) $range(x') - range(x) \leq \ell$.  
\end{lemma}  

\proof 
Wlog, we can restrict ourselves by 
an interval of length 1, 
that is, by a GM-move $x \to x'$. 
Assume for contradiction that 
$range(x) = \ell$, while $range(x') = \ell+1$. 
This happens only if    
% (a) $\max = \max(x_i | i \in [n])$  is kept, while  
% (i) $\min = \min(x_i | i \in [n])$  is reduced by 1. 
$x_n$  is kept unchanged, while $x_1$  is reduced by 1.  
% By (a), $\max$ is a multiple of  $\ell$. 
% Hence, $\min$ is a multiple of  $\ell$ too, 
Hence, $x_n$ is a multiple of $\ell$, since $n$ is the pivot, 
while $x_1$  is a multiple of $\ell$ too, 
since  $range(x) = \ell$. 
But then the pivot is 1, not  $n$, which is a contradiction.  

\smallskip 

Obviously, (ii) immediately follows from (i).  
\qed 

\medskip 
 
In contrast to (ii), $range(x) - range(x')$ 
can be arbitrarily large; see Section~\ref{ss-racing}. 

\subsection{Main theorem. Statement and comments}
% Given a GM-sequence $S = (x^0 \to x^1 \to \ldots \to x^j \to \ldots)$, 
% partition non-negative integers  = \{0,1, \ldots, j, \ldots\}$   
% into an initial part  
% $S_0  = \{0,1, \ldots, j_0\}$  and the infinite part 
% $S_\infty  = \{j_0 + 1, j_0 + 2, \ldots\}$ 
% such that  $range(x^j) \leq \ell$  if and only if  $j \in S_\infty$. 
% By Lemma \ref{l2}, such partition exists.
% In Section \ref{ss-racing} we will show that $j_0$  is bounded 
% by a polynomial in $x$ and $\ell$. 

Let us assign to a GM-sequence $S$  an infinite table 
with rows $j \in \ZZ_+$  and colums $i \in [n]$. 
To each pair $(i,j)$  assign a letter $l = l(i,j)$ 
that can be either $s = s(i,j)$,     
if  $j$  is the pivot in $x^j$,    
or $r = r(i,j)$, if not. 
By definition, $s$  appears exactly once 
in every row  $j$ of $S$. 

Let  $W_i$  be the infinite word 
in the alphabet $\{s,r\}$  corresponding to $i \in [n]$. 
This word is called periodic if there is 
an integer $p > 0$  such that 
$l(i,j) = l(i, j + p)$  for all  $j \geq 0$. 
Clearly, periodicity is respected 
when $p$  is replaced by a multiple of $p$. 
By default we assume that  $p$  is minimal.  
Any interval $P$  of $p$  successive symbols 
will be called a {\em period}. 
We will say that $p = |P|$ is its {\em length}.   
We call two periodic words  $W_i$  and $W_{i'}$  
are equivalent if their periods 
$P_i$  and $P_{i'}$  can be chosen equal.

\begin{theorem} 
\label{t-main} 
For a GM-sequence 
$S = \{x^0 \to x^1 \to \ldots \to x^j \to \ldots\}$ 
such that $range(x^0) \leq \ell$ 
the following statements hold: 

\begin{itemize} 
\item[(1)] 
For every $i \in [n]$ and $j \geq 0$ we have 
$x_i^j -x_i^{j + n \ell} = (n-1) \ell$. 
\item[(2)]  
All infinite words $W_i, i \in [n],$  are equivalent and periodic; 
their common period  $P$  is of length $p = n \ell$;
the corresponding word $W(P)$ 
contains  $\ell$  symbols $s$  and  $(n-1) \ell$  symbols  $r$. 
\item[(3)]  
The length of a maximal sequence of successive $r$'s in $W_i$ 
is a multiple of  $\ell$. 
\item[(4)] 
The length of a maximal sequence of successive $s$' in $W_i$  
is at most $\ell$. 
\end{itemize} 
\end{theorem} 

Obviosly, using (1), one can compute 
$x^j$  for large  $j$ in time 
linear in  $\log(j)$, $\omega(x)$, and $\ell$.   
% \sum_{i \in [n]}(\log(|x_i| + 1)$, 
Furthermore, (1) immediately follows from (2). 

Note that Theorem \ref{t-main} holds even for $n = 1$.  
Indeed, in this case the unique entry $x_1$  is a multiple of $\ell$, 
it remains constant all time, according to the GM-rule.  
Respectivegly,  $\ell (n-1) = 0$. 

\subsection{Left shift lemma. Statement, comments, examples, and proof}
Parts (2,4) of the main theorem 
result easily from the following main lemma.   

\begin{lemma} 
\label{l-main} 
For a GM-sequence 
$S = \{x^0 \to x^1 \to \ldots \to x^j \to \ldots\}$ 
such that $range(x^0) \leq \ell$, we have 
$pivot(x^j) - pivot(x^{j+\ell}) = 1$ for every $j \geq 0$.
\end{lemma} 

This Lemma means that pivots are periodical with period  $\ell$.  
More precisely, the pivot is reduced by 1 after every $\ell$ GM-moves. 
Standardly, we assume that summation is $\mod n$, 
that is, $1 - 1 = 0 = n$.  

Consider $n \ell \times n$ subtable $S'$ of $S$  
formed by the first $n \ell$ rows of $S$  and $i \in [n]$. 
Parition  $S'$ into  $n$  subtables  % $\ell \times n$-subtables  
$S_m$ containing $\ell$ successive rows 
$\{m \ell, \ldots, m \ell +\ell - 1\}$  for    
$m \in \{0,1, \ldots, n-1\}$.   
If we identify $n$  these  $(\ell \times n)$-subtables 
then $n \ell$  pivots partition the obtained 
$(\ell \times n)$-table, by Lemma \ref{l-main}.

Yet, (3) does not result from these observations. 
Not every row-selection in an $\ell \times n$-subtable of $S$ 
can be realized by pivots. 
For example, if several pivots appear in the same 
column  $i \in [n]$, the corresponding rows must be successive, 
since any break contradicts (3).
We will prove (3) in Section \ref{ss-main-proof}.

\medskip 

Let us consider several examples 
with $n=4$ and  $\ell = 3$ (see Table \ref{tab:4subt}).

\begin{table}[h] %htbp
\setlength\tabcolsep{0pt} % let LaTeX figure out intercol. whitespace
\small % 10% linear reduction in font size

\caption{Examples to the Left Shift Lemma for $n=4$ and  $\ell = 3$\\ ("pi" means "pivot")} 
\label{tab:4subt}

\begin{subtable}{0.22\textwidth}
\begin{tabular*}{\linewidth}{@{\extracolsep{\fill}}ccccc}
\toprule
$x_1$ & $x_2$ & $x_3$ & $x_4$  &  pi \\
\midrule
15& 15& 17& 18 & 2  \\
14& 15& 16& 17 & 2  \\
13& 15& 15& 16 & 3  \\
\hline
12& 14& 15& 15 & 1  \\
12& 13& 14& 14 & 1  \\
12& 12& 13& 13 & 2  \\
\hline
11 &12& 12& 12 & 4  \\
10& 11& 11& 12 & 4 \\
9& 10 &10 &12  & 1 \\
\hline
9& 9& 9& 11 & 3  \\
8& 8& 9& 10 & 3 \\
7& 7& 9& 9  & 4 \\
\hline
6& 6& 8 &9 &  2  \\
\bottomrule
\end{tabular*}
\center \scriptsize {$W(P)=(s^2,r^3,s,r^6)$} \label{tab:sub1}
\end{subtable}
\hfill
\begin{subtable}{0.22\textwidth}
\begin{tabular*}{\linewidth}{@{\extracolsep{\fill}}ccccc}
\toprule
$x_1$ & $x_2$ & $x_3$ & $x_4$  &  pi\\
\midrule
15& 16& 17& 17 & 1 \\
15& 15& 16& 16 & 2 \\
14& 15& 15& 15 & 4 \\ 
\hline
13& 14& 14& 15& 4 \\
12& 13& 13& 15& 1\\ 
12& 12& 12& 14& 3\\
\hline
11& 11& 12& 13& 3\\  
10& 10& 12& 13& 4\\
9 & 9& 11& 12& 2\\
\hline
8 &9 &10 &11 & 2\\ 
7 &9 &9 &10 & 3\\
6 &8 &9 &9  &1\\
\hline
6 &7 &8 &8 &1 \\  
\bottomrule
\end{tabular*}
\center \scriptsize {$W(P)=(s^2, r^3, s, r^6)$} \label{tab:sub2}
\end{subtable}%
\hfill
\begin{subtable}{0.22\textwidth}
\begin{tabular*}{\linewidth}{@{\extracolsep{\fill}}ccccc}
\toprule
$x_1$ & $x_2$ & $x_3$ & $x_4$  &  pi \\
\midrule
15& 17& 17& 18&1 \\
15& 16 &16 &17&1 \\
15& 15& 15& 16 & 3 \\
\hline
14& 14& 15 &15& 4 \\
13& 13& 14 &15&  4 \\
12& 12& 13& 15 & 2 \\
\hline
11 &12& 12& 14&3  \\
10& 11& 12& 13 &3\\
9 & 10& 12& 12& 1\\
\hline
9 &9& 11 &11& 2 \\
9 &9 &10& 10 &2\\
7 &9& 9& 9 & 4\\
\hline
6& 8 &8& 9& 1\\
\bottomrule
\end{tabular*}
\center \scriptsize {$W(P) =(s^2, r^6, s, r^3)$} \label{tab:sub3}
\end{subtable}
\hfill
\begin{subtable}{0.22\textwidth}
\begin{tabular*}{\linewidth}{@{\extracolsep{\fill}}ccccc} 
\toprule
$x_1$ & $x_2$ & $x_3$ & $x_4$  &  pi \\
\midrule
15 &  18 &  18  & 18  &  1  \\
15 & 17  & 17 & 17 &  1  \\
15 &  16 &  16  & 16  &  1  \\
\hline
15 &  15  & 15 &  15  &   4  \\
14  & 14  & 14  & 15 &   4  \\
13  & 13  & 13 & 15  &   4  \\
\hline
12  & 12 &  12 &  15 &  3   \\
11  & 11 &  12  & 14 &  3 \\
10 &  10  & 12 &  13   &  3 \\
\hline
9  & 9 &  12 &  12  &   2 \\ 
8  & 9 &  11 &  11  &   2 \\
7  & 9  & 10 &  10  &  2 \\
\hline
6  & 9  & 9 &  9 & 1  \\ 
\bottomrule
\end{tabular*}
\center \scriptsize {$W(P) =(s^3, r^9)$} \label{tab:sub4}
\end{subtable}%
\hfill
\end{table}

\normalsize

\bigskip 

% \subsection*{Proof of the left shift lemma} 
{\em Proof of the left shift lemma.} 
Suppose that  $i$  is the pivot. %  in some row of  $S$. 
By the GM-rule, $x_i = a \ell$ for some integer $a$.  
If  $x_{i'} = a \ell$ then  $i' \leq i$, 
otherwise  $i' > i$  and it would be the pivot, instead of $i$. 
For the same reason, no  $x_{i'}$  equals $(a-1) \ell$.
Thus, $(a-1) \ell < x_{i'} \leq a \ell$  for all $i' \leq i$.

If $i' > i$ then $a \ell < x_{i'} \leq (a+1) \ell$. 
Indeed, $a \ell = x_i \leq x_{i'}$  results from monotonicity. 
If equality holds then $i'$  is a povot rather than $i$.  
Furthermore, $x_{i'} \leq (a+1) \ell$, 
since $range(x) \leq \ell$. 
For the same reasons, if  $x_{i'} = (a+1) \ell$ then  
$x_{i''} = a \ell$  for all $i'' \leq i$.   

Watch  $x_{i-1}$  for the next $\ell$  GM-moves. 
It is reduced by 1 with each move 
until it reaches $(a-1) \ell$, 
since all this time the pivotal entry equals $a \ell$. 

Yet, in at most  $\ell$  GM-moves  
$x_{i-1}$  will be reduced to $(a-1) \ell$, 
since $x_{i-1} \leq x_i = a \ell$.  
Then $i-1$ becomes the pivot and 
remains in this status 
until the considered  $\ell$  GM-moves end. 
This results from the following observations:    

If  $i' > i-1$  then $x_{i'}$  
cannot be reduced to $(a-1) \ell$ in $\ell$ GM-moves. 
Indeed, if  $i' > i$  then  $x_{i'} > a \ell$, as we know. 
If  $i' = i$  then  $x_{i'} = a \ell$, 
but still it cannot be reduced to $(a-1) \ell$
in $\ell$  GM-moves, 
because the first of them ``is lost", 
since  $i' = i$ is the pivot in the beginning. 

If  $i' < i - 1$  then $x_{i'} \geq (a-1) \ell$, 
since $x_i = a \ell$ and  $range(x) \leq \ell$. 
If  $x_{i'} = (a-1) \ell$ 
then $i'$  is the pivot rather than $i$. 
Hence,  $x_{i'} > (a-1) \ell$  and 
it cannot be reduced to 
$(a-2) \ell$  in $\ell$ GM-moves. 
Thus, after $\ell$  GM-moves $i-1$  will be the pivot. 

\medskip 

Formally, the above arguments work for $i=1$ 
but, since this case is ``a bit special", 
we consider it separately. 
Suppose that  $1$  is the pivot. % in some row of $S$.
By the GM-rule,  $x_1 = a \ell$  for some integer $a$.
Let us show that $a \ell < x_i \leq (a+1) \ell$ 
for all $i > 1$.
Indeed, inequality $x_i \geq a \ell$ holds, by monotonicity. 
If equality holds then $i$  is the pivot rather than 1. 
Hence, $x_i > a \ell$. 
Furthermore, $x_i \leq (a+1) \ell$, 
since $range(x) = x_n - x_1  \leq \ell$. 

Watch $x_n$ for the next $\ell$  GM-moves.
It is reduced by 1 with each move  
until $x_n$  reaches $a \ell$, 
since  all this time 
the pivotal entry equals $a \ell$.  
Yet, in at most  $\ell$  GM-moves  
$x_n$  will be reduced to $a \ell$, 
since $x_n \leq (a+1) \ell$.  
Then $n$ becomes the pivot and 
it remains in this status 
until the considered  $\ell$  GM-moves end, 
since no entry can be reduced to $(a-1) \ell$.
Indeed, if  $i > 1$  then  $x_i > a \ell$, as we know, 
and $range(x) \leq \ell$.  
Although  $x_1 = a \ell$, 
but still it cannot be reduced to $(a-1) \ell$
in $\ell$  GM-moves, 
because the first of them ``is lost", 
since  $1$ is the pivot in the beginning. 
\qed 

\subsection{Main theorem. Proof} 
\label{ss-main-proof}
As we know, (2) implies (1).
The main lemma implies (2)  and (4).
Indeed, choose $\ell$ pivots in the first $\ell$ rows, 
$0, 1, \ldots, \ell -1$, 
one for each row, otherwise arbitrarily.  

\begin{remark} 
Of course, there are restrictions.  
For example (3) must hold. 
Yet, we will not characterize 
feasible $\ell$-selections, 
since (2,4) can be proven without it.
\end{remark} 

Main lemma claims that the next $\ell$  pivots, 
in rows $\ell, \ell+1, \ldots, 2\ell-1$,
are obtained by shifting the first $\ell$ by (-1).
% in the cyclic $\mod n$, order, that is, $1 - 1 = 0 = n$. 
Hence, in $n$  such steps   
the original $\ell$-selection of pivots will be obtained.
If we shifted them in the same  $\ell \times n$-table 
these $n$ shifts would partition it.
Thus, $p = n \ell$ and the remaining claims of (2), 
as well as (4) follow. 

Let us prove  (3).
According to the GM-rule, for every  $j \in \ZZ_+$,   
symbol $s$ appears as the $j$-th letter of  $W_i$ 
for exactly one  $i \in [n]$. 

% This happens in two cases: 
% (j)  either there exists a $j \in [n]$  
% such that  $x'_j = b \ell$  and  $b < a$, or
% (jj) there exists a $j \in [n]$  
% such that  $x'_j = a \ell$  and  $i < j$, 
 
\begin{lemma}
\label{l3}
Let  $i$  be the pivot in row $j$ but not in row $j+1$. 
Then,  $i$  becomes the pivot again 
after $c \ell$  GM-moves, that is, in row $j + c\ell$.
\end{lemma} 

Consider the following example with  $n=5$ and $\ell = 7$. 

\medskip 
\noindent
$(5,5,7,8,9) \to  (4,4,7,7,8)  \to  (3,3,6,7,7)  \to (2,2,5,6,7)  \to  (1,1,4,5,7) \to$ 
\newline
$(0,0,3,4,7) \to (-1,0,2,3,6) \to  (-2,0,1,2,5) \to (-3,0,0,1,4).$

\medskip 

% In the first and last vectors $i=3$  is the pivot, 
% in the second one, it is already not, 
% while in the ninth one, that is, 
% in 7 GM-moves  $i = 3$  becomes the pivot again. 

\proof  
Consider two successive moves  
$x^j  \to x^{j+1}$, 
$x^{j+1} \to  x^{j+2}$ in a GM-sequence 
that begins in $x^0$.  
Suppose that  $i \in [n]$   is the pivot in $x^j$  
but not in $x^{j+1}$.
Then, $x^{j+1}_i = a \ell$ , while  $x^{j+2}_i = a \ell - 1$. 
Furthermore  $i$  can become the pivot again 
only when $x_i^{j+j'+1}$  becomes a multiple of  $\ell$, 
that is, $x_i^{j+j'+1} = (a - c) \ell$ with $a > c > 0$.
This requires $c \ell$  moves. 
Thus, $j'$  is a multiple of $\ell$ too. 
\qed 

\subsection{Main theorem. Examples} 
Consider several examples with the same initial vector 
$x^0 = (16,17,20,20,21)$;  $n =range(x^0) = 5$, 
and  $\ell = 2,3,4,5,7$. 
% recall that $\ell \geq 2$.    

\begin{itemize} 
\item
$\ell = 2$, 
the length  $N$  of the initial interval $S_0$ is 4;
the ranges of vectors $x^j, j \in [4]$  are 
$4, 4, 4, 3$, all larger than $\ell = 2$; 
$x^5 = (14, 14, 15, 15, 16)$  
is the first vector of the GM-sequence of range at most 2,  
% $range(x^5) = \ell = 2$;  
all  $x^j, j \geq 5,$  are of ranges 1 or 2, 
that is, at most $\ell$;  
the period is of length  $p = n \ell =  5 \times 2  = 10$; 
respectively, $x^{7+10} = x^{17} = (6, 6, 7, 7, 8)$; 
$x^5 - x^{17} = (n-1) \ell {\bf e} = 8 {\bf e}$, 
where {\bf e}  is the vector of all ones;    
$W(P) = (s, r^4, s, r^4)$.
A power after a letter standardly denotes 
how many times it is repeated.

\item
$\ell = 3$, $N = 6$;  
the ranges of vectors $x^j, j \in [6]$,  are 
$4, 5, 5, 5, 5, 4$,  larger than $\ell = 3$; \;  
$x^7 = (12, 13, 13, 13, 15)$  is the first vector of range at most 3;  
all  $x^j, j \geq 7,$  are of ranges 2 or 3, that is at most $\ell$; \;   
$p = n \ell = 5 \times 3  = 15$; 
respectively, $x^{7+15} = x^{22} = (0, 1, 1, 1, 3)$; \;  
$x^7 - x^{22} = (n-1) \ell {\bf e} = 12 {\bf e}$; \; 
$W(P) = (s^2, r^3, s, r^9)$.  

\item
$\ell = 4$, $N = 1$, 
$range(x^1) = range(x^0) = 5 > 4 = \ell$; \;  
$x^2 = (15, 16, 18, 18, 19)$ is the first vector of range at most 4;  
the ranges of vectors  $x^j, j \geq 2,$  
take values 1, 2, 3, and 4, that is, at most $\ell$; \;   
$p = n \ell = 5 \times 4  = 20$; 
respectively,  $x^{22} = (-1, 0, 2, 2, 3)$; \; 
$x^2 - x^{22} = (n-1) \ell {\bf e} = 16 {\bf e}$; \; 
$W(P) = (s^3, r^8, s, r^8)$.   

\item
$\ell = 5$, $N = 0$, 
the shortest initial interval is empty, since 
$5 = \ell \geq range(x^0) = 5$;  
$x^0$ is already of range at most $\ell$; \;     
all ranges take values 2, 3, 4 and 5; \;   
$p = n \ell = 5 \times 5  = 25$; 
respectively, $x^{25} = (-4, -3, 0, 0, 1)$; \;  
$x^0 - x^{25} = (n-1) \ell {\bf e} = 20 {\bf e}$; \; 
$W(P) = (s^3, r^5, s, r^{10}, s, r^5)$. 

\item
$\ell = 7$,  $N = 0$, 
the shortest initial interval is empty, since 
$\ell \geq range(x^0)$; \;  $x^0$ is already of range at most $ell$; \;     
all ranges take values 1, 2, 3, 4 and 5; \;   
$p = n \ell = 5 \times 7  = 35$; 
respectively, $x^{35} = (-12, -11, -8, -8, -7)$; \;  
$x^0 - x^{35} = (n-1) \ell {\bf e} = 28 {\bf e}$; \; 
$W_(P) = (s^3, r^7, s, r^7, s^3, r^{14})$.  
\end{itemize} 

Note that range 0 is missing in the above examples, 
but it can appear.  
For example, $x^0$  can be of range 0  or 
$(0,1, 1) \to (0,0,0)$.  
Since, $\ell \geq 2$, by definition,    
vectors of ranges  0, 1, or 2  
cannot be in $S_0$, but are always in $S_\infty$ .  

\subsection{Case: $range(x^0) > \ell$. Pursuit of the leaders}
\label{ss-racing} 
If  $range(x^0) > \ell$  then the GM-sequence  
$S$ beginning in  $x^0$  can be partitioned 
into  $S_0$  and  $S_\infty$. 
In other words, there exists $N \geq 0$ such that 
$range(x^0) > \ell$  if $j \leq N$  and 
$range(x^0) \leq \ell$  if $j > N$.  
We provide here an algorithm computing 
$N+1$  and  $x^{N+1}$ 
(the first vector of $S_\infty$)  
in time linear in  $\log N$, $\ell$, and $\omega$. 
% $\sum_{i \in [n]}(\log(|x_i| + 1)$,   
  
Given $x = (x_1, \ldots x_n)$,  
an $i \in [n]$ is called a {\em leader} 
if  $x_i - x_1 \leq \ell$. 
Let $m = m(x, \ell)$  be the maximal leader. 
Then $[m] = \{1, \ldots, m\}$ is the set of leaders. 
Our algorithm is based on 
an iterative expanding this set.

Consider the GM-sequence that begins with $x^0$. 
Obviously, $[m(x^0, \ell)] = [n]$  
if and only if  $range(x^0) \leq \ell$. 
This case was considered in the previous two sections. 
Here we assume that $range(x^0) > \ell$, 
or in other words, the set of leaders in $x^0$ 
% $[m(x^0, \ell)]$ 
is a proper subset of $[n]$.  

It may happen that  $x_i^0$ 
is a multiple of $\ell$ for no $i \in [m]$. 
In this case, none of them is the pivot  
and we need at most $\ell$ GM-moves 
until such an entry appears. 

Let us assume that $x^0$ has it already, 
just to simplify notation.
Then, the pivot will always remain in $[m]$
until entry $m+1$
(perhaps, together with some others) 
catches up with the leaders. 
It will definitely happen 
and we can easily compute how soon. 

Consider the first $m \ell$ moves 
of the GM-sequence that begins in $x^0$. 
By Theorem \ref{t-main},   
the entry $x^0_i$ will be reduced 
by  $(m -1) \ell$ if $i \leq m$,  
and by  $m \ell$ if  $i > m$. 
Thus, outsiders catch up with the leaders 
by  $\ell$  in every $m \ell$ GM-moves. 
Hence, in $m \ell \lfloor \frac{x_{m+1} - x_m}{\ell} \rfloor$  
GM-moves, the distance between them will be at most  $\ell$ 
and  in the next $\ell$ GM-moves, entry 
$m+1$, perhaps with some others, will join the leaders. 
Thus, we obtain a new initial vector and proceed iterating. 
In at most $n - m$  such iterations 
we obtain the last vector $x^N$  
such that $range(x^N) \leq \ell$, 
that is, the last vector of $S_0$. 

\subsection{Finish line}
Applications to impartial games 
require the following generalization. 
Fix an integer  $d$  such that $0 < d \leq n$,   
introduce a finish level, say 0, and assume that 
GM-sequence stops as soon as 
at least  $d$  entries of  $x$  become non-positive. 
Assume that $x^0$ already contains 
$d_0$ such entries and $m_0$ leaders. 
If $d_0 \geq d$ then GM-sequence is finished before it starts. 
Assume also that  $d_0 < d$  and proceed with the 
iterative process described in the previous section. 
We have to decide what will happen first: 

(i) $d$  entries become non-negative 
and the GM-sequence stops, or 

(ii) some new entries will join the leaders and $m$ increases. 

To do so, we compare 
$N_1 = m \ell \lfloor \frac{x_d}{\ell} \rfloor$ and 
$N_2 = m \ell \lfloor \frac{x_{m+1} - x_m}{\ell} \rfloor$. 

(j) If $N_1 \geq N_2$, we compute $x^{N_2}$ 
and proceed further with at most $\ell$  steps to finish 
the GM-sequence; 

(jj)  If $N_1 < N_2$, we expand the set of leaders and 
proceed with obtained  $m' > m$. 

Note that option (ii) disappears in case $range(x^0) \leq \ell$. 

\section{Case $\ell=2$ and remoteness function of exact slow NIM}
\label{s-games}
We assume that the reader is familiar 
with basic concepts of impartial game theory
(see e.g., \cite{ANW07,BCG01-04} for an introduction) 
and also with the recent paper \cite{GMMV23}, 
where the game NIM$(n,k)$ of exact slow NIM,   
was analyzed for the case  $n = k+1$.  
Its mis\`ere version was considered in \cite{GMMN23}.)

\subsection{Exact Slow NIM} 
\label{ss00}
Game Exact Slow NIM  
was introduced in \cite{GH15} as follows: 
Given two integers  $n$ and $k$  
such that  $0 < k \leq n$  and  
$n$  piles containing  $x_1, \ldots, x_n$  stones.  
By one move it is allowed to reduce any  $k$  piles 
by exactly one stone each. 
Two players alternate turns. 
One who has to move but cannot is the loser 
in the normal version of the game and 
the winner of its mis\`ere version. 
In \cite{GH15}, this game was denoted  NIM$^1_=(n,k)$.
Here we will simplify this notation to  NIM$(n,k)$.

Game NIM$(n,k)$  is trivial if  $k = 1$  or  $k = n$. 
In the first case it ends after $x_1 + \ldots + x_n$ moves  
and in the second one---after  $\min(x_1, \ldots, x_n)$ 
moves. In both cases nothing depends on players' skills. 
All other cases are more complicated. 

The game was solved for  $k=2$  and  $n \leq 6$.
In \cite{GHHC20}, an explicit formula 
for the Sprague-Grundy 
(SG) function was found for  $n \leq 4$, 
for both the normal and mis\`ere versions. 
This formula allows us to compute  the SG function in linear time. 
Then, in \cite{CGKPV21} the P-positions 
of the normal version were found for  $n \leq 6$.  
For the subgame with even  $x_1 + \ldots + x_n$ 
a simple formula for the P-positions was obtained, 
which allows to verify in linear time 
if  $x$ is a P-position and, if not, 
to find a move from it to a P-position.  
The subgame with odd  $x_1 + \ldots + x_n$ is more difficult. 
Still a (more sophisticated)  formula 
for the P-positions was found,  
providing a linear time recognition algorithm.

\subsection{Case $n = k+1$}
\label{ss01}
In \cite{GMMV23} the normal version of the game  
was solved in case $n = k+1$  by the following simple rule: 
\begin{itemize}
\item[(o)] if all piles are odd, 
keep a largest one and reduce all other;
\item[(e)] if there exist even piles, 
keep a smallest one of them and reduce all other.
\end{itemize}

This rule was called the {\em M-rule} in \cite{GMMV23}. 
Obviously, it coincides with the GM-rule,  
if we restrict ourselves 
by non-negative positions $x$  
with at least one even entry.

Given  $x$,  assume that both players follow the M-rule 
and  denote by $\M(x)$  the number of moves 
from  $x$  to a terminal position 
(which has at least two entries 0) 
In \cite{GMMV23} it was proven that  $\M = \R$, where 
$\R$  is the {\em remoteness function} 
introduced by Smith in \cite{Smi66}. 
Thus, M-rule solves the game and, moreover, 
it allows to win as fast as possible in an N-position and 
to resist as long as possible in a P-position.

A polynomial algorithm computing  $\M = \R$  
(and in particular, the P-positions) is given, 
even if $n$ is a part of the input and 
integers are presented in binary form.
The results of the present paper,  
restricted to $\ell = 2$, provide much simpler 
and more efficient algorithms.

Let us also note that 
an explicit formula for the P-positions is known 
only for $n \leq 4$  
and already for $n=3$ %% and $n=4$ 
it is pretty complicated 
\cite{GHHC20}\cite[Appendix]{GMMV23}.

\subsection{Related versions of NIM} 
By definition, the present game NIM$(n,k)$ 
is the exact slow version  
of the famous Moore's NIM$_k$ \cite{Moo910}.  
In the latter game a player, by one move, 
reduces arbitrarily 
(not necessarily by one stone) 
at most $k$ piles from $n$. 

The case  $k=1$  corresponds to the classical NIM  
whose P-position was found by Bouton 
\cite{Bou901}  for both the normal and mis\`ere versions. 

\begin{remark} 
Actually, the Sprague-Grundy (SG) values
of NIM were also computed in Bouton's paper, 
although were not defined explicitly in general. This 
was done later by Sprague \cite{Spr36} 
and Grundy \cite{Gru39} for 
arbitrary disjunctive compounds 
of impartial games; see also \cite{Con76,Smi66}. 

In fact, the concept of a P-position 
was also introduced by Bouton in \cite{Bou901}, 
but only for the (acyclic) digraph of NIM, 
not for all impartial games. 
In its turn, this is a special case 
of the concept of a kernel, 
which was introduced for arbitrary digraphs 
by von Neumann and Morgenstern \cite{NM44}. 

Also the mis\`ere version 
was introduced by Bouton in \cite{Bou901},  
but only for NIM, not for all impartial games. 
The latter was done by Grundy and Smith 
\cite{GS56}; see also \cite{Con76,Gur07,Gur07a,GH18,Smi66}.
\end{remark} 

Moore \cite{Moo910} obtained an elegant explicit formula 
for the P-positions of NIM$_k$ 
generalizing the Bouton's case  $k=1$. 
Even more generally, the positions 
of the SG-values 0 and 1  were efficiently characterized by 
Jenkins and Mayberry \cite{JM80}; see also 
\cite[Section 4]{BGHMM17}.  
Also in \cite{JM80};
the SG function of NIM$_k$  
was computed explicitly for the case  $n = k+1$ 
(in addition to the case $k=1$.)  
In general, no explicit formula, 
nor even a polynomial algorithm, 
computing the SG-values 
(larger that 1)  is known. 
The smallest open case: 
2-values for  $n=4$ and $k=2$. 

\medskip 

The remoteness function of  $k$-NIM was recently studied 
in \cite{BGMV23}. 

\medskip 

Let us also mention 
the exact (but not slow) game NIM$^1_=(n,k)$  
\cite{BGHM15,BGHMM17} 
in which exactly $k$ from $n$ piles are reduced 
(by an arbitrary number of stones) in a move. 
The SG-function was efficiently computed 
in \cite{BGHM15,BGHMM17} for $n \leq 2k$. 
Otherwise, even a polynomial algorithm 
looking for the P-positions is not known
(unless $k=1$, of course). 
The smallest open case is  $n=5$  and  $k=2$.

\subsection{Game NIM$(n,n-1)$  for $\ell$ players} 
The exact slow NIM$(n,k)$, with $k = n-1$  
can be played by $\ell$  players  % $j_1, \ldots, j_\ell$  
as follows. 
A position is a non-negative $n$-vector  $x$. 
% $x = (x_1, \ldots, x_n)$. 
Players make moves in a given cyclic order. 
By one move a player can choose an arbirary 
entry and keep it unchanged 
reducicing $n-1$ remaining entries by 1, 
provided they are positive. 
If  $x$  has at least two non-positive entries, 
it is called a {\em terminal} position, 
the game is over, and the player 
who has to move in $x$ (but cannot) lost, 
while  $\ell-1$  others won. 
The payoffs are defined below. 
A sequence of moves 
from the initial position  $x^0$ 
to a terminal one is called a {\em play}. 
Let  $L(P)$  denote the length, 
that is, the number of moves, of a play  $P$. 
Choose a large constant $C$; 
it should be larger than the length 
of any play from  $x^0$.  
Then the loser pays $C-L(P)$ 
and $\ell-1$  winners share this amount, 
genning $\frac{C - L(P)}{\ell-1}$ each. 
The GM-rule defines the unique strategy foreach player.
If a position  $x$  has no entries which are multiple of $\ell$, 
the GM-move in $x$ keeps the lagest entry and reduces 
the smaller $n-1$ by 1 each. 

\begin{conjecture}
\label{con-1} 
The set of $\ell$ GM-strategies form a (uniform) Nash equilibrium. 
\end{conjecture}

By definition, the GM-strategies are uniform, 
that is, independent of  $x^0$. 
Hence, a  Nash equilibrium in GM-strategies  
is uniform, if exists. 

For $\ell = 2$  this conjecture immediately 
follows from the results of \cite{GMMV23}. 

\section{Concluding remarks and open problems} 
GM-rule defines a deterministic dynamic system. 
It is well known that such system can demonstrate 
a ``chaotic behavior". 
However, it is very simple in the considered case, 
due to Lemma \ref{s-main} and Theorem \ref{s-main}.  
After the first $N$, every next $n \ell$  GM-moves  
reduce all entries of $x$  by the same constant $(n-1)\ell$.  
Furthermore, $N$ and $x^N$  also can be efficiently computed.

Given a non-negative $x$, 
integer $\ell \geq 2$ and $d \geq 1$,  
we can determine in linear time, 
how many M-moves is required to get 
$x'$  with at least $d$  non-positive entries. 
In case $\ell=d=2$, this number 
is the value of the remoteness function  $\M(x) = \R(x)$ 
of the slow NIM game NIM$(n,k)$  with $k=n-1$.   
This game can be extended to the case of $\ell \geq 2$ players; 
see Conjecture \ref{con-1}.

Interestingly the obtained optimal GM-strategies 
are uniformly optimal in the following sense. 
Let us replace  0  by an arbitary 
integer $c$, positive or negative, and  
require that at least  $n - d$ entries of  $x$ are at least $c$. 
The game is over as soon as at least 
$d$  entries become at most  $c$. 
In the obtained game 
the GM-strategies are optimal and 
the same for any {\em even} $c$. 
Which other impartial games admit such 
uniformly optimal strategies? 

A more general open question:     
which other discrete dynamic systems 
are related to impartial games;  
in particular, it remains open 
already for the systems defined 
by the GM-rule with $\ell > 2$. 

\bigskip 

\noindent 
{\bf Acknowledgements.}
% This research was prepared within the framework 
% of the HSE University Basic Research Program.  
This research was supported by Russian Science Foundation, 
grant 20-11-20203, https://rscf.ru/en/project/20-11-20203/

\end{document}